\documentclass[10pt,reqno]{article}
\usepackage{amssymb}
\usepackage{amsthm}
\usepackage[tbtags]{amsmath}
\usepackage{doc}
\usepackage{latexsym}
\usepackage{amscd}
\usepackage[frame,cmtip,arrow,matrix,line,graph,curve]{xy}
\usepackage{graphics}
\usepackage{epsfig}
\usepackage{eucal}
\usepackage{mathrsfs}
\usepackage{cite}
\usepackage{textcomp}

\newcommand*{\QEDA}{\hfill\ensuremath{\square}}
\font\headd=cmr8
\begin{document}
\thispagestyle{plain}
 \markboth{}{}
\begin{center}
\noindent{\large\bf Metallic structures on tangent bundle }\\
\end{center}
\begin{center}
\noindent{\bf Mohammad Nazrul Islam Khan}\\ 
{Department of Computer Engineering, College of Computer, Qassim University,\\
P.O.Box 6688, Buraidah 51452,
Saudi Arabia.\\m.nazrul@qu.eu.sa, mnazrul@rediffmail.com \\
orcid.org/0000-0002-9652-0355}
\end{center}
\begin{center}
\noindent
\end{center}
{\bf Abstract.} The purpose of the present work is to study the complete and horizontal lifts of the metallic structure on tangent bundles with respect to almost product structure. We also establish fundamental formulae related to integrability and horizontal lifts of metallic structures on tangent bundles. Moreover, the study revealed the behavior on cross-section of metallic structure in $M$ to tangent bundle $TM$.\\



\pagestyle{myheadings}
 \markboth{\headd .$~~~~~~~~~~~~~~~~~~~~~~~~~~~~~~~~~~~~~~~~~~~~~\,$}
 {\headd $~~~~~~~~~~~~~~~~~~~~~~~~~~~~~~~~~~~~~$Metallic structures }
\noindent{\bf 1. Introduction}
\setcounter{equation}{0}
\renewcommand{\theequation}{1.\arabic{equation}}
\vspace{0.1in}

A general quadratic equation $x^2-\alpha x-\beta I=0,\alpha, \beta$   are positive integers. The set of the positive solutions denoted by $\sigma^\beta_\alpha=\frac{1}{2}(\alpha+\sqrt{\alpha^2+4\beta})$, is called Metallic Means Family (MMF). The member of this family are, for example, the Gold Mean, the Silver Mean, the Bronze Mean, the Copper Mean, the Nickel Mean and many others. The devlopment of the Metallic Means Family are creditted by Kappraff \cite{TGB,BMA,MPO} and Spinadel \cite{TFO,TMM,FTG,TSO}.

For differnt values of $\alpha, \beta$ in $\sigma^\beta_\alpha=\frac{1}{2}(\alpha+\sqrt{\alpha^2+4\beta})$, we have\\ 
  (i)  The values of $\alpha=1, \beta=1$  gives the Fibonacci sequence 1, 1, 2, 3, 5,\ldots and $\sigma^{1}_{1}=\frac{(1+\sqrt{5}}{2})$(Gold Mean), which has been used in a significant of ancient cultures as a proportion basis to compose music, devise sculptures and painting or construct temples and palaces \cite{FTG}.\\
	(ii) The values of $\alpha=2, \beta=1$  gives the Pell sequence 1, 1, 3, 7, 17\ldots and $\sigma^{1}_{2}=1+\sqrt{2}$ (Silver Mean), its importance is implicated in Design, Architecture, and Physics.\\
(iii)  The values of $\alpha=3, \beta=1$  gives $\sigma^{1}_{3}=\frac{(3+\sqrt{13}}{2})$ (Bronze Mean), which plays pivotal role in studying topics such as dynamical systems and quasi-crystals \cite{TMM}.\\
(iv) The values of $\alpha=4, \beta=1$  gives $\sigma^{1}_{4}=2+\sqrt{5}$  (Subtle Mean).\\
(v) The values of $\alpha=1, \beta=2$  gives $\sigma^{2}_{1}=2$ (Copper Mean).\\ 
(vi) The values of $\alpha=1, \beta=3$  gives $\sigma^{3}_{1}=\frac{(1+\sqrt{13}}{2})$ (Nickel Mean), and so on.

On the other hand, Crasmareanu and Hretcanu introduced and studied metallic structures on Riemannian manifold which are polynomial structures satisfying $J^2-pJ-qI$, where $p$ and $q$ are positive integers \cite{MSO}.
Same authors studied the golden structure on manifold by using corresponding almost product structure\cite{GDG}. Recently geometry of these structures has been studied in \cite{OIO, OSI, AOT, PGS}.

The approach of lift has significant role in novel differentiable geometry. By using the lift function, it is convenient to generalized to differentiable structures on any manifold $M$ to its tangent bundle\cite{LOH}. The complete, vertical, horizontal lifts of tensor field and connections on any manifold $M$ to its tangent bundle $TM$ has been obtained by Yano and Ishihara \cite{TAC}. Das and the author \cite{ARS} have been studied almost product structure by means of the complex, vertical and horizontal lifts of an almost r-contact structure.

The paper is structured as follows: In Section 2, definition and properties of the metallic structure are given. Section 3 is devoted to the study of the complete lift of the metallic structure in the tangent bundles and give an example of triple structure i.e. almost product structures, almost tangent structures and almost complex structures. In Section 4, we study integrability of the metallic structure on tangent bundles. We define horizontal lift of the metallic structure in the tangent bundle in Section 5. Finally, we study lifts of metallic structure on cross-section and establish some results related to Nijenhuis tensor on metallic structure.
		
\vspace{0.3in}
\setcounter{equation}{0}
\renewcommand{\theequation}{2.\arabic{equation}}
\noindent{\bf 2. The metallic structure on manifolds}
\vspace{0.1in}

Let $M$ be a differentiable manifold of class $C^\infty$. A tensor field $\Psi$ of type (1, 1) on $M$  is called the metallic structure on $M$ if $\Psi$ satisfies the equation 
\begin{equation}\label{MS}
\Psi^2-\alpha \Psi-\beta I=0 
\end{equation}
where $\alpha, \beta$ are positive integers.

\textbf{Theorem 2.1} \cite{GDG, DGO} An almost product structure $P$ induces the metallic structure as 
\begin{equation}\label{MEP}
\Psi=\frac{1}{2}(\alpha+\sqrt{\alpha^2+4\beta}P)
\end{equation}
Conversely, a metallic structure $\Psi$ gives almost product structure $P$ on $M$ then
\begin{equation}\label{PMS}
P=\frac{2\Psi-\alpha}{\sqrt{\alpha^2+4\beta}}.
\end{equation}\QEDA\\\\																									
Let us define two projection operators $r$ and $s$ as\\
\begin{equation}\label{CP}
r=\frac{1}{2}(I+P),  ~~~~~~ s=\frac{1}{2}(I-P)
\end{equation}
where $P$ is an almost product structure. By writing $\Psi=\frac{1}{2}(I+\sqrt{\alpha^2+4\beta}P)$, we obtain
\begin{equation}\label{NCP}
r=\frac{1}{\sqrt{\alpha^2+4\beta}}\Psi-\frac{\alpha-\sigma^\beta_\alpha}{\sqrt{\alpha^2+4\beta}}I,  ~~~~~~ s=-\frac{1}{\sqrt{\alpha^2+4\beta}}\Psi+\frac{\sigma^\beta_\alpha}{\sqrt{\alpha^2+4\beta}}I
\end{equation}
where $\sigma^\beta_\alpha=\frac{1}{2}(\alpha+\sqrt{\alpha^2+4\beta})$, which is solution of the equation $\Psi^2-\alpha\Psi-\beta I$.
We have  	
\begin{equation}\label{CPR}
r+s=I,  ~~~~   rs=sr=0,  ~~~~   r^2=r, ~~~~    s^2=s.
\end{equation}\\

Let $R$ and $S$ two complementary distributions in $M$ corresponding to projection operators $r$ and $s$, respectively \cite{GDG, DGO}. The projection operators $r$ and $s$ fulfill the accompanying relations:
\begin{equation}\label{MCPR}
\Psi r=r\Psi=\sigma^\beta_\alpha r=\frac{\sigma^\beta_\alpha}{\sqrt{\alpha^2+4\beta}}\Psi-\frac{\beta}{\sqrt{\alpha^2+4\beta}}I 
\end{equation}
\begin{equation}\label{MCPS}
\Psi s=s\Psi=(\alpha-\sigma^\beta_\alpha s)=\frac{\sigma^\beta_\alpha+\alpha}{\sqrt{\alpha^2+4\beta}}\Psi+\frac{\beta}{\sqrt{\alpha^2+4\beta}}I
\end{equation}

\vspace{0.3in}
\setcounter{equation}{0}
\renewcommand{\theequation}{3.\arabic{equation}}
\noindent{\bf 3. Complete lifts of metallic structures in the tangent bundle}
\vspace{0.1in}

In differentiable manifold $M$, $T_p(M)$ is the tangent space at a point $p$ of $M$ i.e. the set of all tangent vectors of $M$ at $p$. the set of all tangent vectors of $M$ at $p$. Then the set $TM=\bigcup_{p\in M}T_p(M)$ is the tangent bundle over the manifold. The complete lift $\Psi^C$ of metallic structure $\Psi$ has the local expression \cite{TAC}
\begin{equation}\label{CLC} 
\Psi^C=
\begin{bmatrix}
{\begin{array}{cc}
\Psi^h_i & 0 \\
\partial\Psi^h_i & \Psi^h_i \\
\end{array} }
\end{bmatrix}
\end{equation} 

\textbf{Proposition 3.1} If $\Psi\in\Im_1^1(M)$ is a metallic structure in $M$, then
its complete lift $\Psi^C$ is also a metallic structure in $TM$.\\
\textit{Proof.} For any $\Psi, G\in\Im^1_1(M)$, we have \cite{TAC} 											
\begin{equation}\label{PGT} 
(\Psi G)^C=\Psi^CG^C
\end{equation} 
If we $G=\Psi$ in (\ref{PGT}), we obtain 
 \begin{equation}\label{CPS} 
(\Psi^2)^C=(\Psi^C)^2
\end{equation} 
Operating the complete lift on (\ref{MS}), we have $(\Psi^2)^C-\alpha \Psi^C-\beta I^C=0 
$. Using (\ref{CPS}) and $I^C=I$, we get 
 \begin{equation}\label{CPAB} 
(\Psi^C)^2-\alpha \Psi^C-\beta I=0
\end{equation} 
This completes the proof.\QEDA\\

 

Operating the complete lift on (\ref{MEP}) and (\ref{PMS}) using corresponding almost product structure, thus we have

\textbf{Theorem 3.2} If $P\in\Im_1^1(M)$ is an almost product structure in $M$. Then its complete lift $P^C$ induces the metallic structures as: 
\begin{equation}\label{APS} 
\Psi^C=\frac{1}{2}(\alpha+\sqrt{\alpha^2+4\beta}P^C)
\end{equation}  
Conversely, If $\Psi\in\Im_1^1(M)$ is a metallic structure. Then its complpete lift $\Psi^C$ yields an almost product structure in $TM$ as 
\begin{equation}
P^C=\frac{2\Psi^C+\alpha}{\sqrt{\alpha^2+4\beta}}.\nonumber
\end{equation}\\
Using equation(\ref{APS}), let us introduce tangent metallic structure and complex metallic structure\cite{DGO}:\\\\ 
\textbf{3.1 Tangent metallic structure}			

Let $T$ be an almost tangent manifold on $M$, then  			
\begin{equation}
\Psi^C_t=\frac{1}{2}(\alpha+\sqrt{\alpha^2+4\beta}T^C)\nonumber
\end{equation}
is said to be tangent metallic structure on $(TM,T^C)$. In particular, $\Psi^C_t$ the polynomial equation $(\Psi^C_t)^2-\alpha\Psi^C_t+\frac{\alpha^2}{4}I$.\\\\
\textbf{3.2 Complex metallic structure}

If $J$ be an almost complex manifold on $M$, then  
\begin{equation}
\Psi^C_j=\frac{1}{2}(\alpha+\sqrt{\alpha^2+4\beta}J^C)\nonumber
\end{equation}
is said to be  the complex metallic structure on $(TM,J^C)$. The tensor field $\Psi^C_j$ satisfies the polynomial equation $(\Psi^C_j)^2-\alpha\Psi^C_j+(\frac{\alpha^2}{4}+\beta)I$.\\

\textbf{Example 3.1} 

	In view of (\ref{MS}), it can be written as \cite{OAB, DGO}
 $$\Psi_{F^C}=\frac{1}{2}(\alpha+\sqrt{\alpha^2+4\beta}F^C), \Psi_{P^C}=\frac{1}{2}(\alpha+\sqrt{\alpha^2+4\beta}P^C), \Psi_{J^C}=\frac{1}{2}(\alpha+\sqrt{\alpha^2+4\beta}J^C)$$ 
where $F,P\in\Im^1_1(M)$ and $J=P\circ F$. We have the relation
$$\sqrt{\alpha^2+4\beta}\Psi_{J^C}=2\Psi_{P^C}\Psi_{F^C}-\alpha\Psi_{P^C}-\alpha\Psi_{F^C}+\alpha \sigma^\beta_\alpha I$$
 																

				
		

\vspace{0.3in}
\setcounter{equation}{0}
\renewcommand{\theequation}{4.\arabic{equation}}
\noindent{\bf 4. Integrability of the metallic structures on tangent bundles}
\vspace{0.1in}

In this section, first we establish the integrability condition of the metallic structure. After that integrability of metallic structures on tangent bundle are studied\\

Let $P, \Psi\in\Im_1^1(M)$. The Nijenhuis tensor  $N_P, N_\Psi\in\Im^2_1(M)$ are defined by \cite{GDG, TAC}
\begin{equation}\label{NP}
N_P(X,Y)=[PX,PY]-P[PX,Y]-P[X,PY]+P^2[PX,PY]
\end{equation}
\begin{equation}\label{NPH}
N_\Psi(X,Y)=[\Psi X,\Psi Y]-\Psi[\Psi X,Y]-\Psi[X,\Psi Y]+\Psi^2[\Psi X,\Psi Y]
\end{equation}
for any $X,Y\in\Im^1_0(M)$, we have the relation \cite{DGO}:
 \begin{equation}\label{NN}
N_P(X,Y)=\frac{4}{\alpha^2+4\beta}N_\Psi(X,Y)
\end{equation}

Suppose that $X,Y\in\Im^1_0(M)$ and $\Psi\in\Im^1_1(M)$, the complete lifts have properties \cite{TAC}": 
\begin{equation}\label{XYC}
(X+Y)^C=X^C+Y^C, [X^C,Y^C]=[X,Y]^C, \Psi^CX^C=(\Psi X)^C.
\end{equation}

From (\ref{NCP}), (\ref{CPR}), (\ref{MCPR}), (\ref{MCPS}), (\ref{PGT}) and (\ref{CPS}), we obtain
\begin{equation}\label{CNC}
r^C=\frac{1}{\sqrt{\alpha^2+4\beta}}\Psi^C-\frac{\alpha-\sigma^\beta_\alpha}{\sqrt{\alpha^2+4\beta}}I, ~~~~  s^C=-\frac{1}{\sqrt{\alpha^2+4\beta}}\Psi^C+\frac{\sigma^\beta_\alpha}{\sqrt{\alpha^2+4\beta}}I
\end{equation}
where $\sigma^\beta_\alpha=\frac{1}{2}(\alpha+\sqrt{\alpha^2+4\beta})$, which is solution of the equation $\Psi^2-\alpha\Psi-\beta I$.
It can easily be seen that   	
\begin{equation}\label{CPR}
r^C+s^C=I,~~~     r^Cs^C=s^Cr^C=0, ~~~    (r^C)^2=r^C,~~~     (s^C)^2=s^C.
\end{equation}
 \begin{equation}\label{MCPR}
\Psi^C r^C =r^C \Psi^C =\sigma^\beta_\alpha r^C , 
\end{equation}
\begin{equation}\label{MCPS}
\Psi^C s^C =s^C \Psi^C =(\alpha-\sigma^\beta_\alpha) s^C
\end{equation}\\
where complete lifts $r^C$ of $r$ and  $s^C$ of $s$
Let $N_P^C$ and $N_\Psi^C$ be Nijenhuis tensors of $P^C$ and $\Psi^C$ in $TM$, respectively. Using equations (\ref{NP}) and (\ref{NPH}) that
\begin{equation}\label{CNP}
N_{P^C}(X^C,Y^C)=[P^CX^C,P^CY^C]-P^C[P^C,Y^C]-P^C[X^C,P^CY^C]+(P^2)^C[P^CX^C,P^CY^C]
\end{equation}
\begin{equation}\label{CNPH}
N_{\Psi^C}(X^C,Y^C)=[\Psi^CX^C,\Psi^CY^C]-\Psi^C[\Psi^CX^C,Y^C]-\Psi^C[X^C,\Psi^CY^C]+(\Psi^2)^C[\Psi^CX^C,\Psi^CY^C]
\end{equation}

\textbf{Proposition 4.1} The complete lift $S^C$ of a distribution $S$ in $TM$ is integrable if and only if $S$ is so in $M$.\\
\textit{Proof.} Let $X$ and $Y$ be vector fields in $M$. The distribution $S$ is integrable iff \cite{GDG}
\begin{equation}\label{RS}
r(sX,sY)=0
\end{equation}
Taking account of equation (\ref{RS}) and using (\ref{XYC}), we find
\begin{equation}\label{CRS}
r^C(s^CX^C,s^CY^C)=0
\end{equation}
where $r^C=(I-s)^C=I-s^C$ . Thus the proposition is proved.
\QEDA\\\\
Hence, we have the following proposition.

\textbf{Proposition 4.2} Let $X$ and $Y$ be vector fields in $M$. If the distribution $S$ be integrable in $M$, i.e. $rN_\Psi(sX,sY)=0$\cite{GDG}. Then $S^C$  is integrable in $TM$ iff
 \begin{equation}\label{CNRS}
r^CN_{\Psi^C}(s^CX^C,s^CY^C)=0
\end{equation} 
\textit{Proof.} Taking account of definition of Nijenhuis tensor $N_{\Psi^C}$  in $TM$ and equation (\ref{NPH}), we have, for any $X,Y\in\Im^1_0(M)$,
 \begin{equation}
N_{\Psi^C}(X^C,Y^C)=[\Psi^CX^C,\Psi^CY^C]-\Psi^C[\Psi^C,Y^C]-\Psi^C[X^C,\Psi^CY^C]+(\Psi^2)^C[\Psi^CX^C,\Psi^CY^C]\nonumber
\end{equation}
Therefore 
\begin{align}\label{CNPHS}
N_{\Psi^C}(s^CX^C,s^CY^C)=[\Psi^Cs^CX^C,\Psi^Cs^CY^C]-\Psi^C[\Psi^Cs^CX^C,s^CY^C]\nonumber\\
-\Psi^C[s^CX^C,\Psi^Cs^CY^C]+(\Psi^2)^C[\Psi^Cs^CX^C,\Psi^Cs^CY^C]
\end{align}
Using the (\ref{CNPHS}) along with (\ref{CPAB}), (\ref{MCPR}) and (\ref{MCPS}), we get
\begin{equation}
N_{\Psi^C}(s^CX^C,s^CY^C)=(-\alpha+2\phi) \Psi^C[s^CX^C,s^CY^C]+(\alpha^2-\alpha\phi+2\beta I)[s^CX^C,s^CY^C]\nonumber
\end{equation}
Multiplying both sides by $\frac{1}{\alpha^2+4\beta}r^C$ and using (\ref{MCPR}) and (\ref{MCPS}), we obtain
\begin{equation}
\frac{1}{\alpha^2+4\beta}r^CN_{\Psi^C}(s^CX^C,s^CY^C)=r^C[s^CX^C,s^CY^C]=(rN_\Psi(sX,sY))^C\nonumber
\end{equation}
Using equation (\ref{CRS}), i.e. $rN_\Psi(sX,sY)=0$ , we have
 \begin{equation}
r^CN_{\Psi^C}(s^CX^C,s^CY^C)=0.\nonumber
\end{equation}
Thus the proposition is proved.\QEDA

\textbf{Proposition 4.3} The complete lift $R^C$ of a distribution $R$ in $TM$ is integrable iff $R$ is integrable in $M$.\\
\textit{Proof.} Let $X$ and $Y$ be vector fields in $M$ and $R$ is integrable if and only if \cite{GDG}
 \begin{equation}\label{SR}
s(rX,rY)=0.
\end{equation}
By taking account of (\ref{SR}) and applying (\ref{XYC}), we have
\begin{equation}\label{CSR}
s^C(r^CX^C,r^CY^C)=0.
\end{equation}
where $s^C=(I-r)^C=I-r^C$. Thus the proposition is proved.\QEDA

 \textbf{Proposition 4.4} Let $X$ and $Y$ be vector fields in $M$ and the distribution $R$ be integrable in $M$,$sN_{\Psi}(rX,rY)=0$\cite{GDG}. Then $R^C$ is integrable in $TM$ iff
 \begin{equation}
s^CN_{\Psi^C}(r^CX^C,r^CY^C)=0.\nonumber
\end{equation}
\textit{Proof.} Let $N_{\Psi^C}$ be the Nijenhuis tensor of $\Psi^C$. Then from (\ref{CNPH}), we have
\begin{align}\label{CNPHR}
N_{\Psi^C}(r^CX^C,r^CY^C)=[\Psi^Cr^CX^C,\Psi^Cr^CY^C]-\Psi^C[\Psi^Cr^CX^C,r^CY^C]\nonumber\\
-\Psi^C[r^CX^C,\Psi^Cr^CY^C]+(\Psi^2)^C[\Psi^Cr^CX^C,\Psi^Cr^CY^C]
\end{align}
The equation (\ref{CNPHR}) along with (\ref{CPAB}), (\ref{MCPR}) and (\ref{MCPS}) gives
 \begin{equation}
N_{\Psi^C}(r^CX^C,r^CY^C)=(\alpha-2\phi) \Psi^C[r^CX^C,r^CY^C]+(\alpha\phi+2\beta I)[r^CX^C,r^CY^C]\nonumber
\end{equation}
Above equation multiplying  by $\frac{1}{\alpha^2+4\beta}s^C$ and applying (\ref{MCPR}) and (\ref{MCPS}), we get
 \begin{equation}
\frac{1}{\alpha^2+4\beta}s^CN_{\Psi^C}(r^CX^C,r^CY^C)=s^C[r^CX^C,r^CY^C]=(sN_\Psi(rX,rY))^C\nonumber
\end{equation}
Using (\ref{CSR}), i.e. $sN_\Psi(rX,rY)$, we have
\begin{equation}
s^CN_{\Psi^C}(r^CX^C,r^CY^C)=0\nonumber
\end{equation}
The proposition is proved.\QEDA

\textbf{Proposition 4.5} Let $X$ and $Y$ be vector fields in $M$ and $\Psi=\frac{1}{2}(\alpha+\sqrt{\alpha^2+4\beta}P)$. Then $N_{P^C}$ and $N_{\Phi^C}$ are related by
 \begin{equation}
N_{P^C}(X^C,Y^C)=\frac{4}{\alpha^2+4\beta}N_{\Psi^C}(X^C,Y^C)\nonumber
\end{equation}
\textit{Proof.} By (\ref{NN}), (\ref{XYC}) and (\ref{CNP}), we have
\begin{equation}
N_{P^C}(X^C,Y^C)=(N_P(X,Y))^C=\left(\frac{4}{\alpha^2+4\beta}N_{\Psi}(X,Y)\right)^C\nonumber
\end{equation}
Thus (\ref{NP}), (\ref{NPH})  and (\ref{XYC}) imply 
 \begin{equation}
N_{P^C}(X^C,Y^C)=\frac{4}{\alpha^2+4\beta}N_{\Psi}^C(X^C,Y^C).\nonumber
\end{equation}
\QEDA

\textbf{Proposition 4.6} Let $\Psi$ be an integrable metallic structure. Then the metallic structure $\Psi^C$ is integrable in $TM$ iff
\begin{equation}
N_{\Psi}^C(X^C,Y^C)=0.\nonumber
\end{equation}
\textit{Proof.} Taking account of definition of  Nijenhuis tensor $N_{\Psi^C}$  in $TM$, we have from (\ref{CNPH})
\begin{equation}
N_{\Psi^C}(X^C,Y^C)=[\Psi^CX^C,\Psi^CY^C]-\Psi^C[\Psi^C,Y^C]-\Psi^C[X^C,\Psi^CY^C]+(\Psi^2)^C[\Psi^CX^C,\Psi^CY^C]\nonumber
\end{equation}
 Using (\ref{XYC}), we obtain
\begin{equation}
N_{\Psi^C}(X^C,Y^C)=(N_\Psi(X,Y))^C=0\nonumber
\end{equation}
since $N_\Psi(X,Y)=0$. 
\QEDA\\
 
\textbf{Example 4.1} 
The complete lifts $R^C$ and $S^C$ of complementary distributions orthogonal $R$ and $S$ corresponding to the Euclidean space of $\Re^2$ are given by
\begin{equation}\label{SPR}
R^C=Span\left\{\frac{\partial}{\partial x}-(x+y)\frac{\partial}{\partial y}\right\}
\end{equation}                                                   
\begin{equation}\label{SPS}
S^C=Span\left\{(x+y)\frac{\partial}{\partial x}+\frac{\partial}{\partial y}\right\}
\end{equation}	
	                                                                   
The equations (\ref{SPR}) and (\ref{SPS}) are associated to the metallic structure		
\begin{equation}\label{ODX}
\Psi\left(\frac{\partial}{\partial x}\right)=\frac{(\alpha-\sigma^\beta_\alpha)(x+y)^2+\sigma^\beta_\alpha}{(x+y)^2+1}\frac{\partial}{\partial x}-\frac{(\sqrt{\alpha^2+4\beta}(x+y)+\sigma^\beta_\alpha}{(x+y)^2+1}\frac{\partial}{\partial y}
\end{equation}                               
\begin{equation}\label{ODY}
\Psi\left(\frac{\partial}{\partial y}\right)=\frac{(\sqrt{\alpha^2+4\beta}(x+y)+\sigma^\beta_\alpha}{(x+y)^2+1}\frac{\partial}{\partial x}+\frac{\sigma^\beta_\alpha(x+y)^2+(\alpha-\sigma^\beta_\alpha)}{(x+y)^2+1}\frac{\partial}{\partial y}
\end{equation}      
which is integrable since $N_\Psi\left(\frac{\partial}{\partial x},\frac{\partial}{\partial y}\right)=0$. 
\newpage
\vspace{0.3in}
\setcounter{equation}{0}
\renewcommand{\theequation}{5.\arabic{equation}}
\noindent{\bf 5.	Horizontal lift of the metallic structure}
\vspace{0.1in}

Suppose that there is a tensor fields $S$ in $M$ and $\nabla_\gamma S$, $\nabla$ is affine connection,  in $TM$ are given by
\begin{equation}
S=S^{i\ldots h}_{k\ldots j}\frac{\partial}{\partial x^i}\otimes\ldots\otimes\frac{\partial}{\partial x^h}\otimes dx^k\otimes\ldots\otimes dx^j\nonumber
\end{equation} 
\begin{equation}
\nabla_\gamma S=y^l\nabla_\gamma S^{i\ldots h}_{k\ldots j}\frac{\partial}{\partial x^i}\otimes\ldots\otimes\frac{\partial}{\partial y^h}\otimes dx^k\otimes\ldots\otimes dx^j\nonumber
\end{equation}
corresponding to the induced coordinates $(x^h,y^h)$ in $\pi^{-1}(U)$\cite{TAC}.\\

Now, we define the horizontal lift $S^H$ of a tensor field $S$ in $M$ to $TM$ by
\begin{equation}
S^H=S^C-\nabla_\gamma S.\nonumber
\end{equation} 	

\textbf{Theorem 5.1} If $\Psi\in\Im^1_0(M)$ is a metallic structure in $M$ then the horizontal lift $\Psi^H$ is also metallic structure in $TM$.\\
\textit{Proof.} For any $\Psi\in\Im^1_1(M)$, we have\cite{TAC}
\begin{equation}\label{PH}
(\Psi^2)^H=(\Psi^H)^2
\end{equation}	 
Operating the horizontal lift on (\ref{MS}), we get $(\Psi^2-\alpha\Psi-\beta I)^H=0$. Taking account of (\ref{PH}) and $I^H=I$, we get
\begin{equation}\label{HMS}
(\Psi^H)^2-\alpha\Psi^H-\beta I=0
\end{equation}	                                                     
we see that $\Psi^H$ is a metallic structure in $TM$. \QEDA
 
Let us introduce a tensor field $\tilde{J}$ of type (1,1) in $TM$ by
\begin{align}\label{JHV}
\tilde{J}X^H=\frac{1}{2}(X^H+\sqrt{\alpha^2+4\beta}X^V)\\
\tilde{J}X^V=\frac{1}{2}(X^V+\sqrt{\alpha^2+4\beta}X^H)
\end{align}	 
then it is easily to show that                                                     
\begin{equation}\label{HMS}
\tilde{J}^2-\alpha\tilde{J}-\beta I=0
\end{equation}	                                                     
Thus $\tilde{J}$ is a metallic structure in $TM$ and we have

\textbf{Theorem 5.2}	
Let $\Psi$ be a metallic structure in $M$ with affine connection$\nabla$. Then $\tilde{J}$ defined by
$$\tilde{J}X^H=\frac{1}{2}(X^H+\sqrt{\alpha^2+4\beta}X^V)$$
$$\tilde{J}X^V=\frac{1}{2}(X^V+\sqrt{\alpha^2+4\beta}X^H)$$ 
is also a metallic structure in $TM$.

\vspace{0.3in}
\setcounter{equation}{0}
\renewcommand{\theequation}{6.\arabic{equation}}
\noindent{\bf 6. Lifts of metallic structure on a cross-section}
\vspace{0.1in}

Let $V$ be a vector field in n-dimensional manifold $M$ and $TM$ its tangent bundle. An n-dimensional submanifold $\beta_V(M)$ of $TM$ is called the cross section determined by $V$, where $\beta_V$ is a mapping $\beta_V:M\rightarrow TM$. If the vector field $V$ has local components $V^h(x)$ in $M$ then the cross section is locally defined by \cite{TAC}
\begin{equation}\label{XYH}
x^h=x^h,  y^h=V^h(x)
\end{equation}	                                                     			
with respect to the induced coordinates $(x^A)=(x^h,y^h)$ in $TM$. 
Let $x^h$ be local component of field $X\in\Im^1_0(M)$ and the local component of vector field $BX$ is 
\begin{equation}\label{BXY}
BX:(B^A_iX^i)=
\begin{bmatrix}
{\begin{array}{cc}
x^h  \\
x^i\partial_iV^h \\
\end{array} }
\end{bmatrix},
\end{equation}	                                                     					
in $TM$, where $BX$ is tangent to $\beta_V(M)$ and defined globally along $\beta_V(M)$. 
The local component of vector field   is 
\begin{equation}\label{CX}
CX:(C^A_iX^i)=
\begin{bmatrix}
{\begin{array}{cc}
0  \\
x^h \\
\end{array} }
\end{bmatrix},
\end{equation}							
which is tangent to the fibre, since a fibre is locally expressed by $x^h=$constant, $y^h=y^h, y^h$ is parameters.\\

From (\ref{BXY}) and (\ref{CX}), we have
\begin{equation}\label{BXBY}
[BX,BY]=B[X,Y], ~~~~~~~~~~~[CX,CY]=0
\end{equation}
for any $X,Y\in\Im^1_0(M)$.	                                                     			
By the definitions of complete and vertical lifts and equations (\ref{BXY}) and (\ref{CX}), we have along $\beta_V(M)$ the formulas
\begin{equation}\label{XBC}
X^C=BX+C(L_VX),~~~~~ X^V=CX
\end{equation}	                                                     
for any $X\in\Im^1_0(M)$, where $L_VX$ denotes the Lie derivative of $X$ with respect to $V$.\\ 
The complete lift $X^C$ and vertical lift $X^V$ of a vector field $X$ in $M$ along $\beta_V(M)$ with local components has components of the form
\begin{equation}\label{XCV}
X^C:
\begin{bmatrix}
{\begin{array}{cc}
x^h  \\
L_Vx^h \\
\end{array} }
\end{bmatrix},
X^V:
\begin{bmatrix}
{\begin{array}{cc}
0  \\
x^h \\
\end{array} }
\end{bmatrix},
\end{equation}\\							
The complete lift $\Psi^C$ of a metallic structure $\Psi$ of $\Im^1_0(M)$ along $\beta_V(M)$ with local components $\Psi^h_i$ in $M$ to $T(M)$ has components of the form
 \begin{equation}\label{PCV}
\Psi^C:
\begin{bmatrix}
{\begin{array}{cc}
\Psi^h_i & 0 \\
L_V\Psi^h_i & \Psi^h_i \\
\end{array} }
\end{bmatrix}
\end{equation}							
then we have along the cross section $\beta_V(M)$ the formula
\begin{equation}\label{PCB}
\Psi^C(BX)=B(\Psi X)+C(L_V\Psi)X
\end{equation}	                                                     
for any $X\in\Im^1_0(M)$. When $\Psi^C(BX)$ is tangent to $\beta_V(M)$, then $\Psi^C$ is said to leave $\beta_V(M)$ invariant. Thus we have

\textbf{Theorem 6.1} \cite{TAC} The complete lift $\Psi^C$  of a metallic structure $\Psi$ of $\Im^1_1(M)$ leave the cross section $\beta_V(M)$ invariant iff $L_V\Psi=0$.

\textbf{Theorem 6.2 } Let $\Psi$ be an almost product structure in $M$ and satisfy the condition $L_V\Psi=0, V$ is a vector field in $M$ then $\Psi^{C\#}$ is a metallic structure on the cross section in $T(M)$ determined by $V$.\\ 
\textit{Proof.} The complete lift $\Psi$ of an element $\Psi$ of  $\Im^1_1(M)$ leave the cross section $\beta_V(M)$ invariant. Let us define an element $F^{C\#}\in\Im^1_1(\beta_V(M))$  by
\begin{equation}\label{PCH}
\Psi^{C\#}(BX)=\Psi^C(BX)=\Psi(BX),   \forall X\in\Im^1_1(\beta_V(M))
\end{equation}	                                                     		
The element $\Psi^{C\#}$ is called the tensor field induced on $\beta_V(M)$ from $\Psi^C.$
since $\Psi$ is a metallic structure in $M$ and $L_V\Psi=0$ i.e 
\begin{equation}\label{PAB}
\Psi^2-\alpha\Psi-\beta I=0~~~~and ~~~~L_V\Psi=0,
\end{equation}	                                                     								
 from (\ref{PCH}), we have
\begin{equation}\label{CPAB}
(\Psi^{C\#})^2-\alpha\Psi^{C\#}-\beta I=0,
\end{equation}	                                                     								
Hence, $\Psi^{C\#}$ is a metallic structure in $\beta_V(M).$
\QEDA\\
Let $N_\Psi$ and $N_{\Psi^C}$ be Nijenhuis tensors of $\Psi\in\Im^1_1(M)$ and of the complete lift $\Psi^C$ of $F$ respectively, from page 36,\cite{TAC}, we obtain
\begin{equation}
N^C_\Psi=(N_\Psi)^C\nonumber
\end{equation}	                                                     		
which implies from (\ref{PCB}), 
\begin{equation}\label{NBC}
N_{\Psi^C}(BX,BY)=B(N_\Psi(X,Y))+C((L_V N_\Psi)(X,Y)
\end{equation}	                                                     
for any $X,Y\in\Im^1_0(M)$. Thus we have from (\ref{NBC})

\textbf{Theorem 6.3}
Let $N_\Psi$ and $N_{\Psi^C}$ be Nijenhuis tensors of $\Psi\in\Im^1_1(M)$ and of the complete lift $\Psi^C$ of $\Psi$ respectively. Then , in order that $N_{\Psi^C}(BX,BY)$ is tangent to the cross-section $\beta_V(M)$ determined by $V\in\Im^1_0(M)$ for any $X,Y\in\Im^1_0(M)$, it is necessary and sufficient that $L_VN_\Psi=0.$\\\\

Suppose that the complete lift $\Psi^C$  of a metallic structure $\Psi$ of $\Im^1_1(M)$ leave the cross section $\beta_V(M)$ invariant. From (\ref{PCH}) and (\ref{BXBY}), we obtain
\begin{align}\label{NPCS}
N_{\Psi^C}(BX,BY)=[\Psi^C(BX),\Psi^C(BY)]-\Psi^C[\Psi^C(BX),(BY)]-\Psi^C[(BX),\Psi^C(BY)]\nonumber\\
+(\Psi^C)^2[\Psi^C(BX),\Psi^C(BY)]\nonumber\\
=[\Psi^{C\#}(BX),\Psi^{C\#}(BY)]-\Psi^{C\#}[\Psi^{C\#}(BX),(BY)]
-\Psi^{C\#}[(BX),\Psi^{C\#}(BY)]\nonumber\\
+(\Psi^{C\#})^2[\Psi^{C\#}(BX),\Psi^{C\#}(BY)]
\end{align}
i.e. 
\begin{equation}\label{NPCSB}
N_{\Psi^C}(BX,BY)
=N_{\Psi^{C\#}}(BX,BY)
\end{equation}
for any $X,Y\in\Im^1_0(M)$, From(\ref{NBC}),
\begin{equation}\label{NBCS}
N_{\Psi^C}(BX,BY)=B(N_\Psi(X,Y))+C((L_V N_\Psi)(X,Y)
\end{equation}	
for any $X,Y\in\Im^1_0(M)$.
As $L_V\Psi=0$ implies that $L_VN_\Psi=0$.
Thus we have\\

\textbf{Theorem 6.4} Let the complete lift $\Psi^C$  of a metallic structure $\Psi$ of $\Im^1_1(M)$ leave the cross section $\beta_V(M)$ invariant. Then 
\begin{equation}
N_{\Psi^{C\#}}=0\nonumber
\end{equation}	
iff
\begin{equation}
N_\Psi=0.\nonumber
\end{equation}


\end{document}